\documentclass[final,authoryear,3p,times]{elsarticle}

\usepackage{amssymb}
\usepackage{lineno}
\usepackage{url}
\usepackage{amsmath}

\journal{Statistics and Probability Letters}

\begin{document}

\begin{frontmatter}

\title{Logarithmic Lambert $\mathrm{W}\times {\cal F}$ random variables for the family of chi-squared distributions and their applications}

\author[um]{Viktor Witkovsk\'y\corref{cor1}}
\ead{witkovsky@gmail.com,witkovsky@savba.sk,umerwitk@savba.sk}
\cortext[cor1]{Corresponding author. Tel.: +421 2 59104530; fax: +421 2 54775943.}
\address[um]{Institute of Measurement Science, Slovak Academy of Sciences, Bratislava, Slovakia}

\author[bratislava,bystrica]{Gejza Wimmer}
\address[bratislava]{Mathematical Institute, Slovak Academy of Sciences, Bratislava, Slovakia}
\address[bystrica]{Faculty of Natural Sciences, Matej Bel University, Bansk\'a Bystrica, Slovakia}
\ead{wimmer@mat.savba.sk}

\author[uk]{Tomy Duby}
\address[uk]{Bicester, Oxfordshire, United Kingdom}
\ead{tomyduby@gmail.com}

\begin{abstract}
We introduce a class of logarithmic Lambert W random variables for a specific family of distributions.  In particular, we characterize the log-Lambert W random variables for chi-squared distributions which naturally appear in the likelihood based inference of normal random variables.  
\end{abstract}

\begin{keyword}
logarithmic Lambert W random variables \sep likelihood based inference \sep exact likelihood ratio test.
\MSC 62F25 \sep 62F30
\end{keyword}

\end{frontmatter}

%\linenumbers

%% Section 1
\section{Introduction}

Originally \cite{Goerg2011, Goerg2012} introduced the Lambert $\mathrm{W}\times {\cal F}$  random variables (RVs) and families of their distributions as a useful tool for modeling skewed and heavy tailed distributions. In particular, for a continuous location--scale family of RVs $X\sim {\cal F}_\vartheta$, parametrized by a vector $\vartheta$, Goerg defined a class of location--scale Lambert $\mathrm{W}\times {\cal F}_\vartheta$ random variables
\begin{equation}\label{eq1}
Y = \left\{U \exp\left( \gamma U\right) \right\}\sigma_X + \mu_X, \quad \gamma \in \mathbb R,
\end{equation}
parametrized by the vector $(\vartheta,\gamma)$, $\mu_X$ and $\sigma_X$ are the location and scale parameters, and  $U = (X-\mu_X)/\sigma_X$. 

The inverse relation to (\ref{eq1}) can be obtained via the multivalued Lambert W function, namely, by the branches of the inverse relation of the function $z = u \exp(u)$, i.e., the Lambert W function satisfies $W(z) \exp(W(z)) = z$, for more details see, e.g., \cite{Corless1996}. 

Here we formally introduce a class of related but different (transformed) RVs and their distributions, the \emph{logarithmic Lambert $\mathrm{W}\times {\cal F}_\vartheta$ RVs} for a specific family of distributions ${\cal F}_\vartheta$ defined on the nonnegative real axis. 

We shall focus on the family of the so-called \emph{log-Lambert W $\times$ chi-squared distributions}, which naturally appear in the statistical likelihood based inference of normal RVs. As we shall illustrate later, a specific type of such RVs plays an important role. In particular, the random variable
\begin{equation}\label{Eynsham}
Y = (Q_\nu - \nu) - \nu \log\left(\frac{Q_\nu}{\nu} \right),
\end{equation}
and its generalization
\begin{equation}\label{Eynsham2}
Y = (Q_\nu - a) - c \log\left(\frac{Q_\nu}{b} \right),
\end{equation}
where $Q_\nu$ is a chi-squared distributed random variable with $\nu$ degrees of freedom, $Q_\nu \sim \chi^2_\nu$, and $a$, $b$, $c$ are further parameters. The random variable (\ref{Eynsham}) will be denoted here as the \emph{standard log-Lambert W $\times\ {\chi^2_\nu}$ random variable}. 

\cite{Stehlik2003}  studied related RVs (and their distributions) to derive the exact distribution of the Kullback-Leibler $I$-divergence in the exponential family with gamma distributed observations. In particular, he derived the cumulative distribution function (CDF) and the probability density function (PDF) of the transformed gamma RVs, which are directly related to the here considered log-Lambert $\mathrm{W}\times {\chi^2_\nu}$ RVs. Stehl{\'\i}k showed that the $I$-divergence from the observed vector to the canonical parameter can be decomposed as a sum of two independent RVs  with known distributions. Since it is related to the likelihood ratio statistics, Stehl{\'\i}k also calculated the exact distribution of the likelihood ratio tests and discussed the optimality of such exact tests. Recently \cite{Stehlik2014} applied the exact distribution of the Kullback-Leibler $I$-divergence in the exponential family for optimal life time testing. 

In this paper we shall present a detailed characterization of the family of the log-Lambert $\mathrm{W}\times {\chi^2_\nu}$ RVs and their distributions together with some applications especially useful for the exact likelihood-based statistical inference. 

We, however, do not consider explicitly the data transformation part, although, as suggested by a reviewer, the analysis and practical importance could benefit from transformation statistics, as it was considered and illustrated in \cite{Goerg2011}, where the emphasis was on the data transformation part. Here we shall focus on applications of the log-Lambert $\mathrm{W}\times {\chi^2_\nu}$ RVs for derivation of some exact LRT (likelihood ratio test) distributions.%

In Section~\ref{sec:2} we shall formally introduce the unifying framework of the logarithmic Lambert W RVs for a specific family of continuous distributions ${\cal F}_\vartheta$,  and further, we shall derive the explicit forms of their CDF and  PDF. Focus of Section~\ref{sec:3} is on the distribution of the log-Lambert $\mathrm{W}\times {\chi^2_\nu}$ RVs and their characterization, including the characteristic function, cumulants, and first moments. We also provide means for computing the distributions of linear combinations of these RVs. Section~\ref{sec:4} illustrates their applicability for the exact (small sample) statistical inference on model parameters based on normally distributed observations.%

%% Section 2
\section{Logarithmic Lambert W random variables for certain families of distributions}\label{sec:2}

Let $X$ be a continuous RV with support on the nonnegative real axis with the probability distribution depending on a (vector) parameter $\vartheta$, i.e.~$X\sim {\cal F}_\vartheta$, where  $\cal F$ indicates a family of distributions. Here we shall consider the following class of transformed RVs,  
\begin{equation}\label{transformation}
Y = g_\theta(X) = \theta_1  -  \theta_2  \log(X) + \theta_3  X,
\end{equation}
where $g_\theta(\cdot)$ is a strictly convex log-linear transformation on $x\geq 0$ for real parameters $\theta = (\theta_1, \theta_2, \theta_3 )$, where $\theta_1 \in{\mathbb R}$, $\theta_2  > 0$, and $\theta_3  >0$. The support of $Y$  is the set $y \in \langle y_{min},\infty)$, where 
\begin{equation}\label{ymin}
y_{min} = g_\theta(x_{y_{min}})  = \theta_1 + \theta_2 -\theta_2\log\left(\frac{\theta_2}{\theta_3}\right),
\end{equation}
with $x_{y_{min}} = \theta_2/\theta_3$, $0< x_{y_{min}}<\infty$. Note that 
\begin{equation}
\frac{Y - \theta_1}{\theta_2} = - \log\left\{ X\exp\left(-\frac{\theta_3}{\theta_2} X\right)\right\}.
\end{equation}
Therefore, the random variable $Y$, defined by (\ref{transformation}), will be called \emph{the log-Lambert W random variable associated with the distribution ${\cal F}_\vartheta$ }(the (minus) log-Lambert $W\times {\cal F}_\vartheta$ RV), and the distribution of $Y$ will be denoted by $Y\sim \mathop{\mathrm{LW}} \left({\cal F}_\vartheta,\theta\right)$.

For illustration, let $X\sim \chi^2_\nu$, then we shall denote the corresponding log-Lambert $\mathrm{W}\times {\chi^2_\nu}$ RV and its distribution by $Y\sim \mathop{\mathrm{LW}}(\chi^2_\nu,\theta)$. In particular, the RV (\ref{Eynsham}) can be expressed in this parametrization as $Y = g_\theta(Q_\nu)$ with $Q_\nu\sim\chi^2_\nu$ and $\theta = \left(\nu(\log(\nu)-1), \nu, 1 \right)$, and consequently with $x_{y_{min}} = \theta_2/\theta_3 = \nu$ and $y_{min} =  \theta_1 + \theta_2 -\theta_2\log\left({\theta_2}/{\theta_3}\right) = \nu(\log(\nu)-1) + \nu - \nu\log\left(\nu\right) = 0$. 

Similarly as in \cite{Goerg2011}, we can define the log-Lambert W RVs $Y\sim \mathop{\mathrm{LW}} \left({\cal F}_\vartheta,\theta\right)$ for other commonly used families of distributions ${\cal F}_\vartheta$ with support on positive real axis. For example, for the gamma and the inverse gamma distribution with the parameters $\alpha$ and $\beta$ we get~$Y\sim \mathop{\mathrm{LW}}\left(\Gamma(\alpha,\beta),\theta\right)$ and $Y\sim \mathop{\mathrm{LW}}\left(\mathrm{inv}\Gamma(\alpha,\beta),\theta\right)$, respectively, and for the Fischer--Snedecor $F$ distribution with $\nu_1$ and $\nu_2$ degrees of freedom, we get $Y\sim \mathop{\mathrm{LW}} \left(F_{\nu_1,\nu_2},\theta\right)$. In all these cases $\theta = (\theta_1,\theta_2,\theta_3)$,  $\theta_1 \in{\mathbb R}$, $\theta_2  > 0$, and $\theta_3  >0$. 

Application of the Lambert W function provides the explicit inverse transformation to (\ref{transformation}). This can be directly used to determine the exact distribution of $Y$, given the distribution of $X$. The cumulative distribution function (CDF) of $Y$ (here denoted by $\mathop{\mathrm{cdf}}\nolimits_Y %\equiv \mathop{\mathrm{cdf}}\nolimits_{\mathop{\mathrm{LW}_{{\cal F}_\vartheta}^\theta}} 
\equiv \mathop{\mathrm{cdf}}\nolimits_{\mathop{\mathrm{LW}}({\cal F}_\vartheta,\theta)}$),  i.e.~%$\mathop{\mathrm{cdf}}\nolimits_{\mathop{\mathrm{LW}_{{\cal F}_\vartheta}^\theta}}(y) 
$\mathop{\mathrm{cdf}}\nolimits_{\mathop{\mathrm{LW}}\left({\cal F}_\vartheta,\theta\right)}(y) = \Pr\left(Y\leq y \,|\, Y\sim \mathop{\mathrm{LW}}\left({\cal F}_\vartheta,\theta\right)\right)$, is given by
\begin{equation}\label{cdfLWF}
	%\mathop{\mathrm{cdf}}\nolimits_{\mathop{\mathrm{LW}_{{\cal F}_\vartheta}^\theta}}
	\mathop{\mathrm{cdf}}\nolimits_{\mathop{\mathrm{LW}}({\cal F}_\vartheta,\theta)}(y) = \mathop{\mathrm{cdf}}\nolimits_{{\cal F}_\vartheta}\left({x}_{U}^{\theta}(y)\right) - \mathop{\mathrm{cdf}}\nolimits_{{\cal F}_\vartheta}\left({x}_{L}^{\theta}(y)\right),
\end{equation}
where $\mathop{\mathrm{cdf}}\nolimits_{{\cal F}_\vartheta}(x) = \Pr\left(X\leq x \,|\, X\sim {\cal F}_\vartheta\right)$ is the CDF of the RV $X$, and ${x}_{L}^{\theta}(y)$ and ${x}_{U}^{\theta}(y)$ are the two distinct real solutions of the equation $y = g_\theta({x})$. In particular, ${x}_{L}^{\theta}(y)$ and ${x}_{U}^{\theta}(y)$ are the solutions on the intervals $(0, x_{y_{min}})$ and $(x_{y_{min}}, \infty)$, respectively, given by 
\begin{eqnarray}\label{xLU}
{x}_{L}^{\theta}(y) = -\frac{\theta_2}{\theta_3} W_0 \left( -\frac{\theta_3}{\theta_2}\exp\left(-\frac{y-\theta_1}{\theta_2} \right)\right), \quad \mathrm{ and }\quad
{x}_{U}^{\theta}(y) = -\frac{\theta_2}{\theta_3} W_{-1} \left( -\frac{\theta_3}{\theta_2}\exp\left(-\frac{y-\theta_1}{\theta_2} \right)\right),
\end{eqnarray}
where $W_0(\cdot)$ and $W_{-1}(\cdot)$ are the two real valued branches of the multivalued \emph{Lambert W function}, i.e.~such function that $z = W(z)\exp(W(z))$, for more detailed discussion see e.g.~\cite{Corless1996} and \cite{Stehlik2003}. Fast numerical implementations of the Lambert W function are available in standard software packages such as \textsc{Matlab}, \textsc{R},  \textsc{Mathematica} or \textsc{Maple}.

Based on the properties of the Lambert W function, see e.g.~Lemma~1 in \cite{Stehlik2003}, note that 
\begin{equation}\label{XLUderivatives}
\frac{d }{d y} {x}_{L}^{\theta}(y)= \frac{{x}_{L}^{\theta}(y)}{\theta_3{x}_{L}^{\theta}(y) -\theta_2}, \quad \mathrm{ and }\quad
\frac{d }{d y} {x}_{U}^{\theta}(y) =\frac{{x}_{U}^{\theta}(y)}{\theta_3{x}_{U}^{\theta}(y) -\theta_2}.
\end{equation}
If $X\sim {\cal F}_\vartheta$ is a continuous RV, then, from (\ref{cdfLWF}) and (\ref{XLUderivatives}), we get that the probability density function (PDF) of $Y\sim \mathop{\mathrm{LW}}\left({\cal F}_\vartheta,\theta\right)$ is given by
\begin{equation}\label{pdfLWF}
\mathop{\mathrm{pdf}}\nolimits_{\mathop{\mathrm{LW}}({\cal F}_\vartheta,\theta)}(y)  = 
\frac{{x}_{L}^{\theta}(y)}{\theta_2- \theta_3{x}_{L}^{\theta}(y)} \mathop{\mathrm{pdf}}\nolimits_{{\cal F}_\vartheta}\left({x}_{L}^{\theta}(y)\right) +
\frac{{x}_{U}^{\theta}(y)}{\theta_3{x}_{U}^{\theta}(y) -\theta_2} \mathop{\mathrm{pdf}}\nolimits_{{\cal F}_\vartheta}\left({x}_{U}^{\theta}(y)\right),
\end{equation}
where $\mathop{\mathrm{pdf}}\nolimits_{{\cal F}_\vartheta}(x)$  denotes the PDF of the RV $X\sim {\cal F}_\vartheta$.

For completeness, by $q_{1-\alpha}$ we shall denote the $(1-\alpha)$-quantile of the distribution $\mathop\mathrm{LW}\left({\cal F}_\vartheta,\theta\right)$, i.e.~such value $q_{1-\alpha}$ that 
$\Pr\left(Y\leq q_{1-\alpha}\,|\,Y\sim \mathop{\mathrm{LW}}\left({\cal F}_\vartheta,\theta\right)\right)= 1-\alpha$, so 
\begin{equation}\label{quantile}
q_{1-\alpha} = \mathop{\mathrm{qf}}\nolimits_{\mathop{\mathrm{LW}}({\cal F}_\vartheta,\theta)}(1-\alpha) 
							 \equiv \mathop{\mathrm{cdf}}\nolimits^{-1}_{\mathop{\mathrm{LW}}({\cal F}_\vartheta,\theta)}(1-\alpha),
\end{equation}
where $\mathop{\mathrm{qf}}\nolimits_{\mathop{\mathrm{LW}}({\cal F}_\vartheta,\theta)}(\cdot)$ denotes the quantile function (QF) of the distribution $\mathop{\mathrm{LW}}\left({\cal F}_\vartheta,\theta\right)$. In general, an analytical solution for $\mathop{\mathrm{qf}}\nolimits_{\mathop{\mathrm{LW}}({\cal F}_\vartheta,\theta)}(\cdot)$ is not available.

% Section 3
\section{Distribution of log-Lambert W random variables for family of chi-squared distributions}\label{sec:3}

Here we  consider the log-Lambert W RV $Y\sim \mathop{\mathrm{LW}}({\chi^2_\nu},\theta)$. Recognizing that the CDF of $\chi^2_\nu$ RV can be expressed by help of incomplete gamma function, directly from (\ref{cdfLWF}) we get
\begin{equation}\label{cdfLWchi2}
	\mathop{\mathrm{cdf}}\nolimits_{\mathop{\mathrm{LW}}({\chi^2_\nu},\theta)}(y) = \frac{1}{\Gamma\left(\frac{\nu}{2}\right)} \Gamma\left(\frac{\nu}{2},{x}_{L}^{\theta}(y),{x}_{U}^{\theta}(y) \right),
\end{equation}
where $\Gamma(\cdot)$ is the \emph{gamma function}, and $\Gamma(n,a,b) = \int_a^b t^{n-1}\exp(-t)\,d t$ is the  \emph{generalized incomplete gamma function}.

Further, as the PDF of the $\chi^2_\nu$ RV is
\begin{equation}\label{pdfchi2}
	\mathop{\mathrm{pdf}}\nolimits_{\chi^2_\nu}(x) = \frac{2^{-\frac{\nu}{2}}}{\Gamma\left(\frac{\nu}{2}\right)} x^{\frac{\nu}{2}-1}\exp\left(-\frac{x}{2}\right),
\end{equation}
directly from (\ref{pdfLWF}) we get
\begin{equation}\label{pdfLWchi2}
\mathop{\mathrm{pdf}}\nolimits_{\mathop{\mathrm{LW}}({\chi^2_\nu},\theta)}(y) = \frac{2^{-\frac{\nu}{2}}}{\Gamma\left(\frac{\nu}{2}\right)} \left(
\frac{\left({x}_{L}^{\theta}(y)\right)^{\frac{\nu}{2}}\exp\left(-\frac{1}{2}{x}_{L}^{\theta}(y) \right)}{\theta_2- \theta_3{x}_{L}^{\theta}(y)} +
\frac{\left({x}_{U}^{\theta}(y)\right)^{\frac{\nu}{2}}\exp\left(-\frac{1}{2}{x}_{U}^{\theta}(y) \right)}{\theta_3{x}_{U}^{\theta}(y)-\theta_2} \right).
\end{equation}

The characteristic function (CF) of $Y\sim \mathop{\mathrm{LW}}({\chi^2_\nu},\theta)$, can be derived by algebraic manipulations directly from its definition, i.e.~$\mathop{\mathrm{cf}}\nolimits_{\mathop{\mathrm{LW}}({\chi^2_\nu},\theta)}(t) = E(\exp(itY)) = E(\exp(itg_\theta(X)))$ where $X\sim\chi^2_\nu$, and is given by
\begin{equation}\label{cfLWchi2}
\mathop{\mathrm{cf}}\nolimits_{\mathop{\mathrm{LW}}({\chi^2_\nu},\theta)}(t) = \frac{2^{-\frac{\nu}{2}}}{\Gamma\left(\frac{\nu}{2}\right)}\frac{\exp(it\theta_1)\Gamma\left(\frac{\nu}{2}-it\theta_2 \right)}{
\left(\frac{1}{2}-it\theta_3 \right)^{\frac{\nu}{2}-it\theta_2 }}.
\end{equation}

The cumulants ($\kappa_j$,  $ \ j = 1,2,\dots$) of $Y\sim \mathop{\mathrm{LW}}({\chi^2_\nu},\theta)$ are readily obtained by expanding the logarithm of the moment generating function (MGF) into a power series. 
Note that
\begin{equation}\label{mgfLWchi2}
\mathop{\mathrm{mgf}}\nolimits_{\mathop{\mathrm{LW}}({\chi^2_\nu},\theta)}(t) = \mathop{\mathrm{cf}}\nolimits_{\mathop{\mathrm{LW}}({\chi^2_\nu},\theta)}(-it)\quad \mathrm{ for }\quad t<\frac{1}{2}\min\left(\frac{\nu}{\theta_2},\frac{1}{\theta_3} \right). 
\end{equation}
Thus
\begin{eqnarray}\label{cumulants}
	\kappa_1 &=& \theta_1 - \theta_2 \log(2) + \nu \theta_3 - \theta_2 \psi^{(0)}\left(\frac{\nu}{2}\right),  \cr
	\kappa_j &=& 2^{j-1} \Gamma(j-1)\theta_3^{j-1}\left(-j \theta_2 + (j-1)\nu \theta_3 \right) + (-1)^j\theta_2^j\psi^{(j-1)}\left( \frac{\nu}{2}\right),
\end{eqnarray}
for $j = 2,3,\dots$, and $\psi^{(m)}(\cdot)$ is the $m$-th order \emph{polygamma function}, i.e.~the $(m+1)$-derivative of the logarithm of the gamma function. 
Hence, the first four basic moment characteristics of this distribution are:
\begin{eqnarray}
\mathop{\mathrm{mean}}     &=& \kappa_1 =  \theta_1 - \theta_2 \log(2) + \nu \theta_3 - \theta_2 \psi\left(\frac{\nu}{2}\right), \\
\mathop{\mathrm{variance}} &=& \kappa_2 = 2 \theta_3\left(-2 \theta_2 + \nu \theta_3 \right) + \theta_2^2\psi^{(1)}\left( \frac{\nu}{2}\right),  \\
\mathop{\mathrm{skewness}} &=&  \frac{\kappa_3}{\kappa_2^{\frac{3}{2}}} =  
\frac{4\theta_3^{2}\left(-3 \theta_2 + 2\nu \theta_3 \right) - \theta_2^3\psi^{(2)}\left( \frac{\nu}{2}\right)}
{\left(2 \theta_3\left(-2 \theta_2 + \nu \theta_3 \right) + \theta_2^2\psi^{(1)}\left( \frac{\nu}{2}\right)\right)^{\frac{3}{2}}},\\
\mathop{\mathrm{kurtosis}} &=& \frac{\kappa_4}{\kappa_2^2} = \frac{16\theta_3^{3}\left(-4 \theta_2 + 3\nu \theta_3 \right) + \theta_2^4\psi^{(3)}\left( \frac{\nu}{2}\right)}{\left(2 \theta_3\left(-2 \theta_2 + \nu \theta_3 \right) + \theta_2^2\psi^{(1)}\left( \frac{\nu}{2}\right)\right)^2}.
	\end{eqnarray}
%

% Subection 3.1
\subsection{The standard log-Lambert $W \times \chi^2_\nu$ distribution}

As mentioned above and as illustrated by examples in the next section, the central role in the likelihood based inference for normally distributed data plays the RV
\begin{equation}\label{LWchi2standard}
Y_\nu = (Q_\nu - \nu) - \nu \log\left(\frac{Q_\nu}{\nu} \right),
\end{equation}
where $Q_\nu \sim \chi^2_\nu$. RV $Y_\nu \sim \mathop{\mathrm{LW}}({\chi^2_\nu},\theta)$, is the special case of RV (\ref{transformation}) with  $\theta = (\nu(\log(\nu)-1),\nu,1)$. We shall call it the  \emph{standard log-Lambert W $\times\ \chi^2_\nu$ RV with $\nu$ degrees of freedom}.

Then, directly from (\ref{cfLWchi2}) we get the characteristic function of the standard log-Lambert $\mathrm{W} \times \chi^2_\nu$ RV as
\begin{equation}\label{cfLWchi2standard}
\mathop{\mathrm{cf}}\nolimits_{Y_\nu}(t) = \frac{2^{-\frac{\nu}{2}}}{\Gamma\left(\frac{\nu}{2}\right)}\frac{\left(\frac{\nu}{e}\right)^{it\nu}\Gamma\left(\frac{\nu}{2}-it\nu \right)}{
\left(\frac{1}{2}-it \right)^{\frac{\nu}{2}-it\nu }},
\end{equation}
and consequently, for $\nu \rightarrow \infty$, we get the convergence of $Y_\nu$  (in distribution) to the chi-squared distribution with 1 degree of freedom, i.e.
\begin{equation}\label{limitDist}
Y_\nu \underset{\nu \rightarrow \infty}{\overset{{\cal D}}{\longrightarrow}} \chi^2_1,
\end{equation}
for more details see~\ref{A1}.
 
As pointed out by one of the reviewers, a specific question for a practitioner is if using the usual rule of thumb, say $n=30$ observations or $\nu = 30$ degrees of freedom, respectively, 
is a good enough approximation for  application of the central limit theorem. Table~\ref{tab1} illustrates how the standard log-Lambert W $\times\ \chi^2_\nu$ distribution, which is an exact null distribution of the likelihood-ratio statistic for testing the hypothesis about the variance parameter based on random sample from normal distribution, differs from the usual (asymptotic) $\chi^2_1$ approximation and also how fast is the convergence to $\chi^2_1$ for $\nu \rightarrow \infty$ (for more details see the Example~1). \cite{Stehlik2003} presents a detailed comparisons with the $\chi^2$-asymptotic of the likelihood-ratio statistic, however in  a different  statistical model with independent observations from the exponential distribution.%

% Table 1
\begin{table}
\caption{The $(1-\alpha)$-quantiles of the standard log-Lambert W $\times\ \chi^2_\nu$ distribution, i.e.~the distribution of the RV (\ref{LWchi2standard}) with $\nu$ degrees of freedom, computed for selected probabilities  $1-\alpha$ and degrees of freedom $\nu$. Note that for $\nu \rightarrow \infty$ we get the chi-squared distribution with 1 degree of freedom.}
\label{tab1}
\begin{center}
\begin{tabular}{lrrrrrrrrr}
\hline

		 				&  \multicolumn{9}{c}{$\nu$}   																																											\\ 
		 						
		 					 \cline{2-10}
		 				
$1-\alpha$  &          1 &          2 &          3 &          5 &         10 &         20 &         30 &        100 &  $\infty$ \\

\hline 

																																																																\\

    0.7000 &     1.4145 &     1.2543 &     1.1951 &     1.1468 &     1.1103 &     1.0922 &     1.0862 &     1.0778 &     1.0742 \\

    0.7500 &     1.7308 &     1.5426 &     1.4713 &     1.4124 &     1.3677 &     1.3454 &     1.3380 &     1.3277 &     1.3233 \\

    0.8000 &     2.1306 &     1.9105 &     1.8245 &     1.7524 &     1.6974 &     1.6698 &     1.6607 &     1.6479 &     1.6424 \\

    0.8500 &     2.6605 &     2.4039 &     2.2993 &     2.2102 &     2.1415 &     2.1069 &     2.0953 &     2.0792 &     2.0723 \\

    0.9000 &     3.4254 &     3.1259 &     2.9968 &     2.8840 &     2.7956 &     2.7506 &     2.7356 &     2.7146 &     2.7055 \\

    0.9500 &     4.7606 &     4.4077 &     4.2418 &     4.0906 &     3.9683 &     3.9053 &     3.8841 &     3.8543 &     3.8415 \\

    0.9750 &     6.1137 &     5.7256 &     5.5301 &     5.3438 &     5.1885 &     5.1070 &     5.0795 &     5.0406 &     5.0239 \\

    0.9900 &     7.9162 &     7.4984 &     7.2734 &     7.0470 &     6.8499 &     6.7441 &     6.7081 &     6.6570 &     6.6349 \\

    0.9990 &    12.4771 &    12.0220 &    11.7549 &    11.4566 &    11.1683 &    11.0035 &    10.9459 &    10.8635 &    10.8276 \\

    0.9999 &    17.0579 &    16.5840 &    16.2977 &    15.9575 &    15.5983 &    15.3792 &    15.3007 &    15.1868 &    15.1367 \\

\hline
\end{tabular}
\end{center}
\end{table}%

% Subection 3.2
\subsection{Computing the distributions of linear combinations of independent log-Lambert $\mathrm{W} \times \chi^2_\nu$ random variables}

The CDF, PDF and QF of the log-Lambert $W \times \chi^2_\nu$ distribution can be numerically evaluated directly from (\ref{cdfLWF}), (\ref{pdfLWF}), and (\ref{quantile}), by using suitable implementation of the Lambert W function. 
 
Numerical evaluation of the distribution of a linear combination of independent log-Lambert $\mathrm{W} \times \chi^2_\nu$ RVs is based on methods similar to those discussed in \cite{Witkovsky2001a} and \cite{Witkovsky2004}, and is closely related to the method for computing the distribution of a linear combination of independent chi-squared RVs suggested by \cite{Imhof1961}, see also \cite{Davies1980}, and also related to the method for computing the distribution of a linear combination of independent inverted gamma variables suggested by \cite{Witkovsky2001b}. The procedure is based on the numerical inversion of the characteristic function, for more details see \cite{GilPelaez1951}. 

Consider thus the random variable $Y=\sum_{j=1}^k \lambda_j Y_j$, a linear combination of independent log-Lambert W RVs $Y_j \sim \mathop{\mathrm{LW}}({\chi^2_{\nu_j}},\theta_j)$ with $\nu_j$, 
degrees of freedom, parameters $\theta_j = (\theta_{j_1}, \theta_{j_2}, \theta_{j_3})$, and real coefficients $\lambda_j$, $j=1,\dots,k$. Let $\mathop{\mathrm{cf}}\nolimits_{Y_j }(t)$ denote the
characteristic function of $Y_j$. The characteristic function of $Y$ is
\begin{equation}\label{charT}
\mathop{\mathrm{cf}}\nolimits_{Y}(t) = \mathop{\mathrm{cf}}\nolimits_{Y_1 }(\lambda_1t)\cdots \mathop{\mathrm{cf}}\nolimits_{Y_k }(\lambda_kt),
\end{equation}
where
\begin{equation}\label{chart}
\mathop{\mathrm{cf}}\nolimits_{Y_j }(\lambda_jt)=
\frac{2^{-\frac{\nu}{2}}}{\Gamma\left(\frac{\nu}{2}\right)}\frac{\exp(i\lambda_jt\theta_{j_1})\Gamma\left(\frac{\nu}{2}-i\lambda_jt\theta_{j_2} \right)}{
\left(\frac{1}{2}-i\lambda_jt\theta_{j_3} \right)^{\frac{\nu}{2}-i\lambda_jt\theta_{j_2} }}.
\end{equation}

The distribution function (CDF) of $Y$, $\mathop{\mathrm{cdf}}\nolimits_{Y}(y) = \Pr\{Y\leq y\}$, is according to the inversion formula due to \cite{GilPelaez1951} given by
\begin{equation}\label{cdfY}
\mathop{\mathrm{cdf}}\nolimits_{Y}(y)  =  \frac{1}{2}-\frac{1}{\pi}\int_0^\infty \Im\left(\frac{e^{-ity}\mathop{\mathrm{cf}}\nolimits_{Y}(t) }{t} \right)\,dt,
\end{equation}
and the PDF is given by
\begin{equation}\label{pdfY}
\mathop{\mathrm{pdf}}\nolimits_{Y}(y) = \frac{1}{\pi}\int_0^\infty \Re\left(e^{-ity}\mathop{\mathrm{cf}}\nolimits_{Y}(t)  \right)\,dt.
\end{equation}

This approach can also be applied to compute the distributions of more general linear combinations of independent RVs, e.g.~with $\chi^2_{\nu_i}$ and $\mathop{\mathrm{LW}}({\chi^2_{\nu_j}},\theta_j)$ distributions. 

The \textsc{Matlab} implementation of the algorithms for computing CDF, PDF, CF, cumulants and QF of the log-Lambert $\mathrm{W} \times \chi^2_\nu$ RVs (resp.~their linear combinations) is currently available at \url{http://www.mathworks.com/matlabcentral/fileexchange/46754-lambertwchi2},  the \textsc{Matlab} Central File Exchange. In future, these algorithms will become a part of a more general \textsc{Matlab} suite of programs (under development) to calculate the tail probabilities (including CDF, PDF, and QF)  of a linear combination of RVs in one of the following classes: (1) class of symmetric RVs containing normal, Student's $t$, uniform and triangular distributions, and (2) class of RVs with support on positive real axis, e.g., the chi-squared and inverse gamma distributions, see \url{http://sourceforge.net/projects/tailprobabilitycalculator/}.

% Section 4
\section{Examples}\label{sec:4}

For illustration, here we present simple examples of the likelihood based inference for normally distributed data, where the distribution of the likelihood ratio test statistic under the null hypothesis can be expressed using the log-Lambert $\mathrm{W} \times \chi^2_\nu$ RVs.

% Subsection 4.1
\subsection{Example 1: Distribution of the LRT statistic for testing a single variance component}

Let $S^2$ be the estimator of the variance parameter $\sigma^2$ (e.g.~the restricted maximum likelihood estimator (REML) of $\sigma^2$, based on a random sample from normally distributed data, $\tilde{Y} \sim N(0,\sigma^2 I)$), such that $\frac{\nu S^2}{\sigma^2}\sim \chi^2_\nu$. The PDF of the RV $S^2$ can be directly derived from the PDF of the chi-squared distribution with $\nu$ degrees of freedom (\ref{pdfchi2}). So the log-likelihood function is 
\begin{equation}\label{loglikS2}
\mathop{\mathrm{loglik}}\left(\sigma^2 \,|\, S^2\right) =  \mathop{\mathrm{const}} + \left(\frac{\nu}{2} - 1\right)\log\left(\frac{\nu S^2}{\sigma^2} \right) - \frac{1}{2}\left( \frac{\nu S^2}{\sigma^2} \right) + \log\left( \frac{\nu}{\sigma^2} \right),
\end{equation}
and the (log-) likelihood-ratio test statistic (LRT) for testing $H_0: \sigma^2 = \sigma^2_0$ vs.~alternative $H_A: \sigma^2 \neq \sigma^2_0$ is 
\begin{equation}\label{LRT1}
\mathop{\mathrm{lrt}} =  -2 \left(\sup_{H_0}\left\{\mathop{\mathrm{loglik}}\left(\sigma^2 \,|\, S^2\right) \right\} - \sup\left\{\mathop{\mathrm{loglik}}\left(\sigma^2 \,|\, S^2\right) \right\}\right) 
= -2 \left( \mathop{\mathrm{loglik}}\left(\sigma^2_0 \,|\, S^2\right) -  \mathop{\mathrm{loglik}}\left(\hat{\sigma}^2 \,|\, S^2\right) \right),
\end{equation}
where $\hat{\sigma}^2 = S^2$ is the REML estimator of $\sigma^2$. From that,
\begin{eqnarray}\label{LRT1b}
\mathop{\mathrm{lrt}}  &=& 
 -2 \left(\left[  \left(\frac{\nu}{2} - 1\right)\log\left(\frac{\nu S^2}{\sigma^2_0} \right) - \frac{1}{2}\left( \frac{\nu S^2}{\sigma^2_0} \right) + \log\left( \frac{\nu}{\sigma^2_0} \right)\right]  -  \left[  \left(\frac{\nu}{2} - 1\right)\log\left(\frac{\nu S^2}{\hat{\sigma}^2} \right) - \frac{1}{2}\left( \frac{\nu S^2}{\hat{\sigma}^2} \right) + \log\left( \frac{\nu}{\hat{\sigma}^2} \right)\right]\right)\cr
&=& \left(\frac{\nu S^2}{\sigma^2_0} - \nu\right) - \nu \log \left(\frac 1\nu \frac{\nu S^2}{\sigma^2_0} \right) \stackrel{H_0}{\sim} (Q_\nu - \nu) - \nu \log\left(\frac{Q_\nu}{\nu} \right),
\end{eqnarray}
where $Q_\nu \sim \chi^2_\nu$. That is, under the null hypothesis $H_0$, the LRT statistic has the standard log-Lambert $\mathrm{W}\times \chi^2_\nu$ distribution,  $\mathop{\mathrm{lrt}}  \stackrel{H_0}{\sim} \mathop{\mathrm{LW}}(\chi^2_\nu,\theta)$ with $\theta = \left(\nu(\log(\nu)-1),\nu,1\right)$. 

Based on that, the $(1-\alpha)$-confidence interval for the parameter $\sigma^2$, say $\mathop{{c}}_{\mathrm{LRT}}$, obtained by inverting the LRT, can be expressed as
\begin{equation}
\mathop{{c}}\nolimits_{\mathrm{LRT}} = \left\{\sigma^2: \left(\frac{\nu S^2}{\sigma^2} - \nu\right) - \nu \log \left(\frac 1\nu\frac{\nu S^2}{\sigma^2} \right)\leq q_{1-\alpha} \right\}
=	\left\{\sigma^2: \frac{\nu S^2}{{x}_{U}^{\theta}(q_{1-\alpha})} \leq \sigma^2 \leq  \frac{\nu S^2}{{x}_{L}^{\theta}(q_{1-\alpha})}\right\},
\end{equation}
where $q_{1-\alpha}$ denotes the $(1-\alpha)$-quantile of the random variable $(Q_\nu - \nu) - \nu \log\left(\frac{Q_\nu}{\nu} \right)$, i.e.~$\mathop{\mathrm{LW}}({\chi^2_\nu},\theta)$ distribution with $\theta = \left(\nu(\log(\nu)-1),\nu,1\right)$, and the limits ${x}_{L}^{\theta}(q_{1-\alpha})$ and ${x}_{U}^{\theta}(q_{1-\alpha})$ are defined by (\ref{xLU}). The minimum length confidence interval for $\sigma^2$, say $\mathop{{c}}\nolimits_{\mathrm{ML}}$,  can be expressed, e.g., as
\begin{equation}
\mathop{{c}}\nolimits_{\mathrm{ML}} = \left\{\sigma^2: \left(\frac{\nu S^2}{\sigma^2} - \nu\right) - (\nu + 2) \log \left(\frac 1\nu\frac{\nu S^2}{\sigma^2} \right)\leq \tilde{q}_{1-\alpha} \right\}
=	\left\{\sigma^2: \frac{\nu S^2}{{x}_{U}^{\tilde{\theta}}(\tilde{q}_{1-\alpha})} \leq \sigma^2 \leq  \frac{\nu S^2}{{x}_{L}^{\tilde{\theta}}(\tilde{q}_{1-\alpha})}\right\},
\end{equation}
where $\tilde{q}_{1-\alpha}$ is the $(1-\alpha)$-quantile of the random variable  $(Q_\nu - \nu) - (\nu+2) \log\left(\frac{Q_\nu}{\nu} \right)$, i.e.~$\mathop{\mathrm{LW}}(\chi^2_\nu,\tilde\theta)$ distribution with $\tilde{\theta} = \left((\nu+2)\log(\nu)-\nu,\nu+2,1\right)$, see also \cite{Tate1969} and \cite{Juola1993}.

% Sunsection 4.2
\subsection{Example 2: Distribution of the LRT statistic for testing normal linear regression model parameters}

Let $Y\sim N(X\beta,\sigma^2 I)$ be an $n$-dimensional normally distributed random vector with a non-stochastic full-ranked $(n\times k)$-design matrix $X$, parameters  $\beta\in{\mathbb R}^k$ and $\sigma^2>0$. Here the log-likelihood function is 
\begin{equation}\label{loglikSigmaML}
\mathop{\mathrm{loglik}}\left(\beta, \sigma^2 \,|\, Y\right) = - \frac{n}{2} \log(2\pi) - \frac{n}{2} \log\left(\sigma^2\right) - \frac{1}{2\sigma^2} (Y - X\beta)^T(Y - X\beta),
\end{equation}
and the LRT statistic for testing $H_0: (\beta,\sigma^2) = (\beta_0,\sigma^2_0)$, vs.~alternative $H_A: (\beta,\sigma^2) \neq (\beta_0,\sigma^2_0)$,  is given by
\begin{eqnarray}\label{LRT2}
\mathop{\mathrm{lrt}} 
&=&  -2 \left( \mathop{\mathrm{loglik}}\left(\beta_0,\sigma^2_0 \,|\, Y\right) -  \mathop{\mathrm{loglik}}\left(\hat{\beta},\hat{\sigma}^2 \,|\, Y\right) \right)
=  \frac{1}{\sigma^2_0} (Y - X\beta_0)^T(Y - X\beta_0) - n\log\left(\frac{\hat{\sigma}^2}{\sigma^2_0} \right)-n \cr
&\stackrel{H_0}{\sim}&  Q_k + \left\{\left(Q_{\nu} - n\right) - n \log\left(\frac{Q_{\nu}}{n}\right)\right\},
\end{eqnarray}
where $\hat{\beta} = (X^T X)^{-1}X^T Y$ and $\hat{\sigma}^2 = \frac{1}{n}(Y - X\hat{\beta})^T(Y - X\hat{\beta})$, such that $Q_\nu = \frac{n \hat{\sigma}^2}{\sigma^2_0} \sim \chi^2_\nu $, with $\nu = n-k$, and independent of $Q_k\sim \chi^2_k$. That is, under the null hypothesis $H_0$, the LRT statistic (\ref{LRT2}) is distributed as a linear combination (sum) of two independent RVs with $\chi^2_k$ and $\mathop{\mathrm{LW}}(\chi^2_\nu,\theta)$ distributions, respectively, where $\nu = n-k$ and $\theta = \left(n(\log(n)-1),n,1\right)$. For more details, see e.g.~\cite{Choudhari2001} and \cite{Chvostekova2009}.

% Subsection 4.3
\subsection{Example 3: Distribution of the (restricted) LRT statistic for testing canonical variance components}

Consider a normal linear model with two variance components, $Y \sim N(X\beta, \sigma_1^2 V + \sigma^2_2 I_n)$, where $Y$ is an $n$-dimensional normally distributed random vector, $X$ is a known $(n\times k)$-design matrix, $\beta$ is a $k$-dimensional unknown  vector of fixed effects, $V$ is a known $n\times n$ positive semi-definite matrix, $I_n$ is the $n\times n$ identity matrix, and $\sigma_1^2\geq 0$, $\sigma^2_2 >0$ are the variance components --- the parameters of interest. 

The (restricted) LRT methods are based on distribution of the maximal invariant $\tilde{Y} = B^TY$, where $B$ is an arbitrary matrix, such that $BB^T %= M_X 
= I_n - X X^+$ (here $X^+$ denotes the Moore-Penrose $g$-inverse of $X$) and $B^TB  = I_{n-\mathop{rank}(X)}= I_\nu$.  Hence, $\tilde{Y} \sim N(0,\Sigma)$, where $\Sigma = \sigma_1^2 W + \sigma^2_2 I_\nu = \sum_{i=1}^r (\sigma_1^2 \varrho_i + \sigma^2_2) D_i = \sum_{i=1}^r \vartheta_i D_i$, $W = B^TVB = \sum_{i=1}^r  \varrho_i D_i$ is a spectral decomposition with the eigenvalues $\varrho_i$ ($\varrho_1 > \cdots > \varrho_r \geq 0$) and their multiplicities $\nu_i$,  $\nu = \sum_{i=1}^r \nu_i$,  $D_i$ are mutually orthogonal symmetric matrices such that $D_iD_i = D_i$, $D_iD_j = 0$ for $i\neq j$, and $I_\nu = \sum_{i=1}^r D_i$. Here we consider the problem of testing hypothesis about canonical variance components $\vartheta = (\vartheta_1,\dots,\vartheta_r)$, where $\vartheta_i = \sigma_1^2 \varrho_i + \sigma^2_2$, $i = 1,\dots,r$. Namely, we consider testing the hypothesis $H_0: \vartheta = \vartheta_0$, vs.~alternative $H_A: \vartheta \neq \vartheta_0$. 

Let $U_i = \tilde{Y}^T D_i \tilde{Y}$, according to \cite{Ohlsen1976}, the following holds true: $U = (U_1,\dots,U_r)$ is a minimal sufficient statistic for the parameters $(\sigma_1^2,\sigma^2_2)$, and $U_i/(\sigma_1^2 \varrho_i + \sigma^2_2) = U_i/\vartheta_i\equiv Q_{\nu_i} \sim \chi^2_{\nu_i}$, $i = 1,\dots,r$ are mutually independent chi-squared RVs with $\nu_i$ degrees of freedom. 
Thus, for specific values of the canonical parameter $\vartheta_0 = ({\vartheta_0}_1,\dots,{\vartheta_0}_r)$ and the minimal sufficient statistic $U = (U_1,\dots,U_r)$, the (restricted) log-likelihood function can be expressed as 
\begin{equation}
\mathop{\mathrm{loglik}}\left(\vartheta_0 \,|\, U\right) = -\frac{\nu}{2}\log(2\pi) - \frac{1}{2}\sum_{i=1}^r\nu_i\log({\vartheta_0}_i) - \frac{1}{2}\sum_{i=1}^r\frac{U_i}{{\vartheta_0}_i}, 
\end{equation}
and the (restricted) LRT statistic is given by
\begin{eqnarray}\label{LRT3}
\mathop{\mathrm{lrt}}  &=& -2 \left(\mathop{\mathrm{loglik}}\left(\vartheta_0 \,|\, U\right) - \mathop{\mathrm{loglik}}\left(\hat{\vartheta} \,|\, U\right)\right)
= \sum_{i=1}^r \left\{\left(\frac{U_i}{{\vartheta_0}_i} - \nu_i\right) - \nu_i \log\left(\frac{1}{\nu_i}\frac{U_i}{{\vartheta_0}_i}\right)\right\} \cr
&\sim& \sum_{i=1}^r \left\{\left(\frac{\vartheta_i}{{\vartheta_0}_i} Q_{\nu_i} - \nu_i\right) - \nu_i \log\left(\frac{\vartheta_i}{{\vartheta_0}_i}\frac{ Q_{\nu_i}}{\nu_i}\right)\right\}  
\stackrel{H_0}{\sim}\sum_{i=1}^r \left\{\left(Q_{\nu_i} - \nu_i\right) - \nu_i \log\left(\frac{Q_{\nu_i}}{\nu_i}\right)\right\}, 
\end{eqnarray}
where $\hat{\vartheta} = (\hat{\vartheta}_1,\dots,\hat{\vartheta}_r)$ with  $\hat{\vartheta}_i = \frac{U_i}{\nu_i}$ is the REML estimator of $\vartheta$,  ${\vartheta} = ({\vartheta}_1,\dots,{\vartheta}_r)$ represents the true (unknown) vector of parameters,  and $Q_{\nu_i}\sim\chi^2_{\nu_i}$ are mutually independent RVs, $i = 1,\dots,r$. 

That is, under the null hypothesis $H_0: \vartheta = \vartheta_0$, the restricted LRT statistic (\ref{LRT3}) is distributed as a linear combination of $r$ independent RVs with $\mathop{\mathrm{LW}}(\chi^2_{\nu_i},{\theta_i})$ distributions, where $\theta_i = \left(\nu_i(\log(\nu_i)-1),\nu_i,1\right)$. 
In general, if $\vartheta \neq \vartheta_0$,  the LRT statistic (\ref{LRT3}) is distributed as a linear combination of $r$ independent RVs with $\mathop{\mathrm{LW}}(\chi^2_{\nu_i},{\theta}_i)$ distributions, where $\theta_i = \left(\nu_i(\log(\nu_i/\lambda_i)-1),\nu_i,\lambda_i\right)$ with $\lambda_i = \frac{\vartheta_i}{{\vartheta_0}_i} $, $i = 1,\dots,r$. %

% Subsection 4.4
\subsection{Example 4: Numerical example} 
% Table 2
\begin{table}
\caption{The (different) eigenvalues $\varrho_i$ of the $W$ matrix from model (\ref{ANOVAmodel}), with their multiplicities $\nu_i$, together with the observed values of the sufficient statistics $U_i$, the true values of the canonical variance components $\vartheta^*_i = \sigma^2_1 \varrho_i + \sigma_2^2$, and the hypothetical values of the parameters $\vartheta^{H_{01}}_i = \sigma^2_1\varrho_i+\sigma^2_2$, $\vartheta^{H_{02}}_i = 0\times\varrho_i+\sigma^2_2$, and $\vartheta^{H_{03}}_i = 1\times\varrho_i+\sigma^2_2$, where $i = 1,\dots,r$, $\sigma^2_1 = 0.1$ and $\sigma^2_1 = 1$.}
\label{tab2}
\begin{center}
\begin{tabular}{rrrrrrrr}
\hline

$\rho_i$    &  $\nu_i$  &     $U_i$ & $\vartheta^*_i$ & $\hat{\vartheta}_i$ & $\vartheta^{H_{01}}_i$ & $\vartheta^{H_{02}}_i$& $\vartheta^{H_{03}}_i$ \\

\hline \\

     19.24 &       1 &       0.65 &       2.92 &       0.65 &       2.92 &       1.00 &      20.24 \\

     17.04 &       1 &      17.12 &       2.70 &      17.12 &       2.70 &       1.00 &      18.04 \\

     14.89 &       1 &       2.76 &       2.49 &       2.76 &       2.49 &       1.00 &      15.89 \\

     12.77 &       1 &       3.01 &       2.28 &       3.01 &       2.28 &       1.00 &      13.77 \\

     10.65 &       1 &       0.45 &       2.06 &       0.45 &       2.06 &       1.00 &      11.65 \\

      8.53 &       1 &       4.02 &       1.85 &       4.02 &       1.85 &       1.00 &       9.53 \\

      6.42 &       1 &       0.52 &       1.64 &       0.52 &       1.64 &       1.00 &       7.42 \\

      4.30 &       1 &       2.06 &       1.43 &       2.06 &       1.43 &       1.00 &       5.30 \\

      2.16 &       1 &       0.90 &       1.22 &       0.90 &       1.22 &       1.00 &       3.16 \\

      0.00 &     100 &     117.25 &       1.00 &       1.17 &       1.00 &       1.00 &       1.00 \\
      
\hline 
\end{tabular}  
\end{center}
\end{table}%

In order to illustrate some of the numerical calculations required for testing hypothesis on canonical variance components based on the LRT statistic, as presented in Example~3, let us consider the following unbalanced one-way random effects ANOVA model, as a special case of a  normal linear model with two variance components:
\begin{equation}\label{ANOVAmodel}
	Y_{ij} = \mu + b_i + \varepsilon_{ij}, \quad i = 1,\dots,G,\  j = 1,\dots,n_i,
\end{equation}
where $Y = (Y_{11},\dots,Y_{Gn_G})^T$ is the $n$-dimensional vector of measurements, $n = \sum_{i=1}^{G} n_i$, $\mu$ represents the common mean, $b = (b_1,\dots,b_G)^T$ is a vector of random effects, $b\sim N(0,\sigma^2_1 I_G)$, and  $\varepsilon = (\varepsilon_{11},\dots,\varepsilon_{Gn_G})^T$ is the $n$-dimensional vector of measurement errors,  $\varepsilon \sim N(0,\sigma^2_2 I_n)$.

%%%%%%%%%%%%%%%%%%%%%%
%Correct !!
% $W = B^T I_n B = \sum_{i=1}^r  \varrho_i D_i$
% Should be
% $W = B^T V B = \sum_{i=1}^r  \varrho_i D_i$
%%%%%%%%%%%%%%%%%%%%%%
In particular, for $G = 10$, and $n_1 = 2, n_2 = 4, n_3 = 6, \dots, n_{10} = 20$, with $n = 110$, by spectral decomposition of the matrix  $W = B^T V B = \sum_{i=1}^r  \varrho_i D_i$, we get $r = 10$  different eigenvalues $\varrho_i$ with their multiplicities $\nu_i$, $i = 1,\dots,r$. For the true values of the parameters $\mu = 0$, $\sigma^2_1 = 0.1$, and $\sigma^2_2 = 1$, we have generated the $n$-dimensional vector of observations $Y$ with the observed values of the sufficient statistics $U = (U_1,\dots,U_r)$. The true values of the canonical variance components, $\vartheta^* = (\vartheta^*_1,\dots,\vartheta^*_r)$ with $\vartheta^*_i = \sigma^2_1 \varrho_i + \sigma_2^2$, were estimated by REML, $\hat{\vartheta} = (\hat{\vartheta}_1,\dots,\hat{\vartheta}_r)$,  where $\hat{\vartheta}_i = U_i/\nu_i$. 

Here the goal is to test the following null hypotheses: $H_{01}: \vartheta = \vartheta^{H_{01}}$ with $\vartheta^{H_{01}}_i = \vartheta^*_i = \sigma^2_1\varrho_i+\sigma^2_2$, as well as  $H_{02}: \vartheta = \vartheta^{H_{02}}$ with $\vartheta^{H_{02}}_i = 0\times\varrho_i+\sigma^2_2$, and $H_{03}: \vartheta = \vartheta^{H_{03}}$ with $\vartheta^{H_{03}}_i = 1\times\varrho_i+\sigma^2_2$, for numerical values see Table~\ref{tab2}.

\begin{figure}
\begin{center}
\includegraphics[width = 0.66\textwidth]{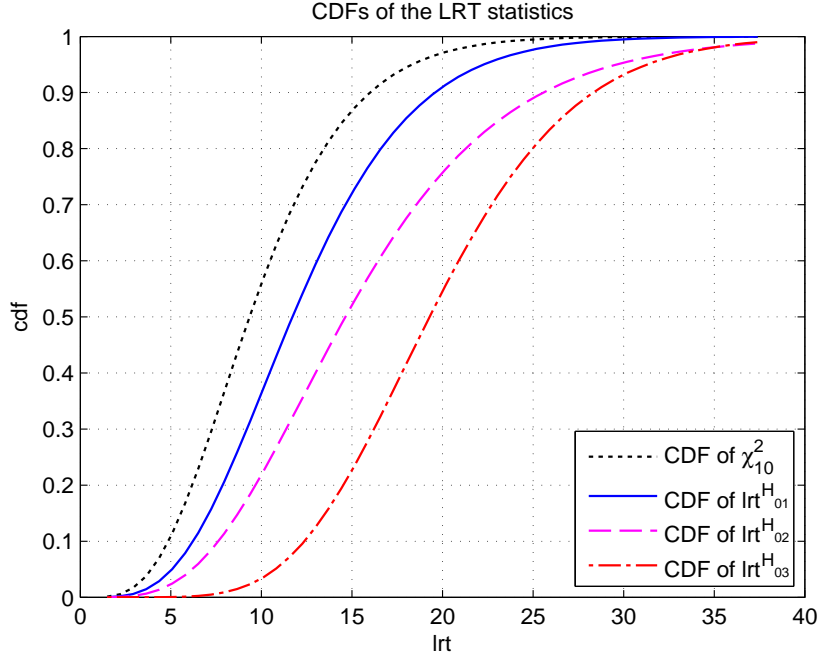}% color figure
\end{center}
\caption{The cumulative distribution functions (CDFs) of the LRT statistics (\ref{LRT3}) for testing the null hypotheses $H_{01}$, $H_{02}$, and $H_{03}$, respectively. 
Here, the true null distribution is the CDF of $\mathop{\mathrm{lrt}}^{H_{01}}$. The CDF of the $\chi^2_{10}$ distribution is frequently used as the approximate (asymptotic) null distribution.}
\label{fig1}%
\end{figure}

The CDFs of the LRT statistics for testing the null hypotheses $H_{01}$, $H_{02}$, and $H_{03}$, respectively, that is $\mathop{\mathrm{lrt}}^{H_{01}}$,  $\mathop{\mathrm{lrt}}^{H_{02}}$, and $\mathop{\mathrm{lrt}}^{H_{03}}$, are plotted in Figure~\ref{fig1}, together with the CDF of the $\chi^2_{10}$ distribution, which is conventionally used as the  approximate (asymptotic) distribution. Note that only $H_{01}$ is true, and so, only $\mathop{\mathrm{lrt}}^{H_{01}}$ has the correct null distribution given by (\ref{LRT3}). 

For given observed values of the sufficient statistics, $U_i$, the observed value of the LRT statistic was $\mathop{\mathrm{lrt}}^{H_{01}} = 7.3095$ ($\mathop{\mathrm{lrt}}^{H_{02}}= 18.7350$, and  $\mathop{\mathrm{lrt}}^{H_{03}}  = 10.6475$). The $(1-\alpha)$-quantile of the null distribution, for $\alpha = 0.05$, is $\mathop{\mathrm{q}}^{H_0}_{1-0.05} = 22.2689$. For comparison, the quantile of the $\chi^2_{10}$ distribution is $\chi^2_{10,1-0.05} = 18.3070$. Based on that, on significance level $\alpha = 0.05$, we cannot reject any of the hypotheses $H_{01}$, $H_{02}$, $H_{03}$. However, note that the hypothesis $H_{02}$ would be rejected if the approximate $\chi^2_{10}$ null distribution were used instead of the exact null distribution.

% Section 5
\section{Conclusions}

In this paper we introduce the class of the log-Lambert $\mathrm{W}\times {\cal F}$ random variables and their distributions. It includes, as special case, the class of log-Lambert $\mathrm{W}\times \chi^2_\nu$ RVs, which naturally appears in statistical inference based on likelihood of normal RVs. A suite of \textsc{Matlab} programs (implementation of algorithms for computing PDF, CDF, QF, CF, cumulants, and convolutions) is available at  \textsc{Matlab} Central File Exchange, \url{http://www.mathworks.com/matlabcentral/fileexchange/46754-lambertwchi2}.

% Section Acknowledgements
\section{Acknowledgements} % \label{}
The work was supported by the Slovak Research and Development Agency, grant APVV-0096-10, SK-AT-0025-12, and by the Scientific Grant Agency of the Ministry of Education of the Slovak Republic and the Slovak Academy of Sciences, grant VEGA 2/0038/12 and VEGA 2/0043/13.


% Section Appendix
\appendix

\section{Limit distribution of the standard log-Lambert $\mathrm{W} \times \chi^2_\nu$ random variables}\label{A1}

To show that (\ref{limitDist}) holds true, note that the moment generating function (MGF) of the standard log-Lambert $\mathrm{W} \times \chi^2_\nu$ RV with $\nu$ degrees of freedom  is
\begin{equation}\label{tag2}
\mathop{\mathrm{mgf}}\nolimits_{Y_\nu}(t)  = \frac{\Gamma \left ( \frac {\nu}{2}-\nu t\right )}{\Gamma \left ( \frac {\nu}{2}\right )} \left ( \frac {\nu}{2e}\right )^{\nu t} (1 - 2t)^{\nu t - \frac {\nu}{2}}, \quad 0\leq t < \frac 12,
\end{equation}
and, according to equation (1.4.24) in \cite{Lebedev1963}, for $x \in \mathbb C$, $|x|\gg 1$, $|\arg x|\leq \frac {\pi}{2}$, we have
\begin{equation}
\Gamma (x) = \sqrt{2 \pi}x^{x-\frac 12} e^{-x}(1+ r(x)),
\end{equation}
where $|r(x)|\leq\frac {c}{|x|}$ for some real $c>0$. Thus, for $\nu \rightarrow \infty$ and for all $t\in \langle 0,\frac 12)$, we get
\begin{eqnarray}
\mathop{\mathrm{mgf}}\nolimits_{Y_\nu}(t) &=& \frac {\Gamma \left ( \frac {\nu}{2}-\nu t\right )}{\Gamma \left ( \frac {\nu}{2}\right )} \left ( \frac {\nu}{2e}\right )^{\nu t} (1 - 2t)^{\nu t - \frac {\nu}{2}}\cr
&=& \frac {\sqrt{2 \pi}\left ( \frac {\nu}{2}-\nu t\right )^{\frac {\nu}{2}-\nu t -\frac 12} e^{- \frac {\nu}{2}+\nu t}\nu^{\nu t} (1 - 2t)^{- \frac {\nu}{2}+\nu t}(1+ r(\frac {\nu}{2}-\nu t))}{\sqrt{2 \pi}\left ( \frac {\nu}{2} \right )^{\frac {\nu}{2}-\frac 12} e^{- \frac {\nu}{2}}2^{\nu t} e^{\nu t}(1+ r(\frac {\nu}{2}))}\cr
&=& \frac {\left ( \frac {\nu}{2}-\nu t\right )^{\frac {\nu}{2}-\nu t -\frac 12} \nu^{\nu t} (1 - 2t)^{- \frac {\nu}{2}+\nu t}(1+ r(\frac {\nu}{2}-\nu t))}{\left ( \frac {\nu}{2} \right )^{\frac {\nu}{2}-\frac 12} 2^{\nu t} (1+ r(\frac {\nu}{2}))}\cr
&=& \left [ \frac {\nu \frac 12 (1-2t)}{\nu\frac 12} \right ]^{\frac {\nu}{2} - \frac 12} \left ( \frac {\nu}{2}-\nu t\right )^{-\nu t } \left (
\frac {\nu}{2}\right )^{\nu t} (1 - 2t)^{- \frac {\nu}{2}+\nu t}\frac {(1+ r(\frac {\nu}{2}-\nu t))}{(1+ r(\frac {\nu}{2}))}\cr
&=& (1 - 2t)^{\frac {\nu}{2}- \frac 12}\left ( \frac {\nu}{2}\right )^{-\nu t } (1 - 2t)^{-\nu t} \left ( \frac {\nu}{2}\right )^{\nu t }                             (1 - 2t)^{\nu t - \frac {\nu}{2}}\frac {(1+ r(\frac {\nu}{2}-\nu t))}{(1+ r(\frac {\nu}{2}))}\cr
&=& \frac {1}{\sqrt {1 - 2t}} \frac {(1+ r(\frac {\nu}{2}-\nu t))}{(1+ r(\frac {\nu}{2}))}.
\end{eqnarray}
Consequently, 
\begin{equation}
\lim_{\nu \to \infty} \mathop{\mathrm{mgf}}\nolimits_{Y_\nu}(t) = \frac {1}{\sqrt {1 - 2t}}, \quad 0\leq t < \frac 12,
\end{equation}
what coincides with the MGF of a chi-squared distribution with 1 degree of freedom.

\end{document}